\documentclass[11pt,a4paper]{article} 
  
\newtheorem{theo}{Theorem}[section]

\newtheorem{cor}[theo]{Corollary} 
\newtheorem{ex}[theo]{Example} 
\newtheorem{defi}[theo]{Definition}

\newcommand{\ot}{\otimes} 
\newcommand{\co}{{\cal C}} 

\def\bra{\langle}
\def\ket{\rangle}

\def\cop{\Delta}

\def\eps{\varepsilon}

\def\du1{\hat 1}

\def\-1{_{(-1)}}
\def\0{_{(0)}}
\def\1{_{(1)}}
\def\2{_{(2)}}
\def\3{_{(3)}}

\def\|{\, | \,}

\def\du1{\hat 1}
\def\lact{\triangleright}

\newcommand{\End}{\mbox{\rm End}\,}

\def\id{\mbox{\rm id}}

\def\to{\rightarrow}

\title{An action-free characterization of weak Hopf-Galois extensions} 
\author{Lars Kadison} 

\begin{document} 
\date{} 
\maketitle 

\begin{abstract} 
We define comodule algebras 
and Galois extensions for  actions of bialgebroids.
Using just module conditions we characterize the Frobenius extensions that are Galois as depth two and right balanced
extensions. As a corollary,
we obtain characterizations of certain weak and ordinary Hopf-Galois extensions
without reference to action in the hypothesis.   
\end{abstract} 



\section{Introduction} 
A finite Hopf $H$-Galois extension $A | B$ is a Frobenius extension
and has many characterizations as an $H^*$-module algebra $A$ with invariants $B$ 
satisfying various conditions, or dually, as an $H$-comodule algebra $A$
with coinvariants $B$ satisfying various conditions \cite{Mo}; and has
many interesting applications \cite{CK,Mo}.
There is recently in \cite{KN,KS} a characterization of
certain noncommutative Hopf-Galois extensions  - those
with trivial centralizer (and arising in subfactor theory) - in terms of module-theoretic
conditions of depth two on the tensor-square of the
 extension and ``balanced'' on the module $A_B$. 
The main thesis in \cite{KS} is that
to a depth two ring extension $A | B$ one associates by construction 
two bialgebroids over the unrestricted centralizer $R$,
a left bialgebroid $S := \End {}_BA_B$ and a right bialgebroid
$T := (A \ot_B A)^B$, the $B$-central elements with multiplication
induced from $T \cong \End {}_A(A \ot_B A)_A$.  The bialgebroids
$S$ and $T$ are simultaneously each other's left and right $R$-dual
bialgebroids and they act on $A$ and $\End {}_BA$ respectively.
If $B$ is trivial, we obtain the two main examples of Lu bialgebroids.
If $R$ is trivial, $R$-bialgebroids are usual bialgebras,
and if $R$ is a separable algebra, $R$-bialgebroids are weak
bialgebras:  antipodes may be added to create weak Hopf
algebras if $A | B$ is additionally a Frobenius
extension. If $A_B$ is balanced, the invariants $A^S = B$
and the endomorphism ring $\End A_B$ is a smash product of $A$ and $S$
(tensoring over $R$), which signals a Galois extension.  Good definitions of Galois
extension in terms of applications
have appeared very recently \cite{BW,CG}.  In this paper
 we extend the main theorems in \cite{KN, KS}
to  weak Hopf-Galois extension and Galois extensions
for bialgebroids.  

\section{Preliminaries} 

Let $B \subseteq A$ be an associative not necessarily commutative subring pair sharing $1$,
also referred to here as a ring extension $A | B$.  
\begin{defi} 
A ring extension $A | B$ is depth two (D2) if
the tensor-square $A \ot_B A$ is isomorphic both as natural $B$-$A$-bimodules
(left D2) and as
$A$-$B$-bimodules (right D2) to a direct summand
of a finite direct sum of $A$ with itself. Equivalently, $A | B$ is D2 if there
are (left D2 quasibase) elements $\beta_i \in S$, $t_i \in T$ such that
\begin{equation}
\label{eq: left}
 a \ot_B a' = \sum_i t_i \beta_i(a)a' 
\end{equation}
and (right D2 quasibase) elements $\gamma_j \in S$, $u_j \in T$ such that 
\begin{equation}
\label{eq: right}
 a \ot_B a' = \sum_j a \gamma_j(a') u_j 
\end{equation}
for all $a,a' \in A$.
\end{defi}

\begin{ex} 
A finite dimensional algebra is D2 with dual
bases as a vector space corresponding to D2 quasibases. 
Given a subgroup of a finite group $H < G$, the complex group subalgebra pair
$\co H \subseteq \co G$ is D2 iff $\, H \, \triangleleft \, G$ \cite{KK}.  
Another related example:
a normal Hopf subalgebra pair is D2. Yet another is  a finite weak
Hopf-Galois extension \cite[3.1]{K2004}.   
\end{ex} 

Recall from \cite{KS} that a right $R'$-bialgebroid $T'$ are two rings $R'$ and $T'$ with two maps
$\tilde{s}, \tilde{t}: R' \to T'$, a ring homomorphism and anti-homomorphism resp.,
such that $\tilde{s}(r) \tilde{t}(r') = \tilde{t}(r') \tilde{s}(r)$ for all $r,r' \in R'$,
$(T',\ \cop \!: T' \to \ $ $ T' \ot_{R'} T', \eps: T' \to R')$ is an $R'$-coring w.r.t. the $R'$-$R'$-bimodule $r \cdot x \cdot r' = x \tilde{t}(r) \tilde{s}(r')$ such that
 $(\tilde{s}(r) \ot 1) \cop(x) = $ $ (1 \ot \tilde{t}(r))\cop(x)$,
$\cop(xy) = \cop(x)\cop(y)$ (which makes sense thanks to the previous axiom), $\cop(1) = 1 \ot 1$, $\eps(1_{T'}) = 1_{R'}$ and
 $\eps(xy) = \eps(\tilde{s}(\eps(x)) y ) = \eps(\tilde{t}(\eps(x)) y)$  for all $x, y \in T', r,r' \in R'$. 
A left bialgebroid is just a right bialgebroid with three of the axioms transposed \cite{KS}. 

\begin{ex} \cite[section 5]{KS}
Given a D2 extension $A | B$, the ring $T =\ $ $ (A \ot_B A)^B$ is a right bialgebroid
over the centralizer $R := C_A(B)$ with maps $\tilde{s}(r) =$ $ 1 \ot r$, $\tilde{t}(r) = r \ot 1$,
$\cop(x) = \sum_j (x^1 \ot_B \gamma_j(x^2)) \ot_R u_j$ and $\eps(x) = x^1x^2$ for all $x = x^1 \ot x^2 \in T$.
Note then that ${}_RT_R$ is given by $r \cdot x \cdot r' = rx^1 \ot x^2 r'$ and ${}_RT$ is finite projective
from eq.~(\ref{eq: right}).  
\end{ex} 
The dual left $R$-bialgebroid is $S$, there being two $R$-valued nondegenerate pairings
of $S$ and $T$; e.g., $\bra \alpha, t \ket = t^1 \alpha(t^2)$ for each $\alpha \in S$, $t \in T$.
The left bialgebroid structure is given by $\overline{s}(r) = \lambda_r$, left multiplication
by $r \in R$, $\overline{t}(r) = \rho_r$, right multiplication, $\cop(\alpha) = \sum_i \alpha(?t_i^1)t_i^2 \ot \beta_i$
and $\eps_S(\alpha) = \alpha(1)$.  

 The $R$-bialgebroid $S$ acts on $A$ by 
$\alpha \lact a = \alpha(a)$ (where $R$ acts as a subring of $A$) with invariants
$A^S = \{ a \in A | \alpha \lact a =$ $\eps_S(\alpha)a \, \forall \alpha \in S \} $ $ \supseteq$ $B$,
and if $A_B$ is balanced, $A^S = B$.  $A$ is thereby a left $S$-module algebra (or algebroid \cite[2.1]{KS}).  We
need the dual notion:  
 
\begin{defi}
Let $T'$ be a right $R'$-bialgebroid $(T', \tilde{s}, \tilde{t},$ $ \cop, $ $ \eps)$.   
A (right) $T'$-comodule algebra $A'$ is a ring $A'$ with ring homomorphism $R' \to A'$ 
 together with a coaction $\delta: A' \to A' \ot_{R'} T'$, where  values $\delta(a)$ are denoted by the Sweedler
notation $a\0 \ot a\1$, 
such that $A'$ is a right $T'$-comodule over the $R'$-coring $T'$ \cite[18.1]{BW}, 
$ \delta(1_{A'}) = 1_{A'} \ot 1_{T'} $, $ra\0 \ot a\1 =$ $ a\0 \ot \tilde{t}(r)a\1$
for all $r \in R'$,
and $\delta(aa') = \delta(a) \delta(a')$ for all $a,a' \in A'$.   The subring of coinvariants
is ${A'}^{\rm co \, T'} := \{ a \in A' | \delta(a) = a \ot 1_{T'} \}$. Consequently $R'$ and ${A'}^{\rm co \, T'}$ commute in $A'$. 
\end{defi} 

For example, $T'$ is a comodule algebra over itself. A D2 extension $A | B$ has $T$-comodule algebra
$A$ \cite[5.1]{K2004}, indeed a $T$-Galois extension, which we define as follows.  

\begin{defi}
Let $T'$ be a left finite projective right $R'$-bialgebroid.  
A  $T'$-comodule algebra $A'$ is a (right) $T'$-Galois extension of its coinvariants
$B'$ if the (Galois) mapping 
$\beta: A' \ot_{B'} A' \to A' \ot_{R'} T'$ defined
by $\beta(a \ot a') = a{a'}\0 \ot {a'}\1 $ is bijective.  
\end{defi}

\section{D2 characterization of Galois extensions}
In this section we provide characterizations of generalized Hopf-Galois extensions
in analogy with the Steinitz characterization of Galois extension of fields as being
separable and normal.  
 
\begin{theo}
\label{th-main}
Let $A \| B$ be a Frobenius extension. The extension $A | B$
is $T$-Galois for some left finite projective right bialgebroid $T$
over some ring $R$ if and only if
$A | B$ is D2 with $A_B$ balanced.   
\end{theo}
\textbf{Proof}. 
($\Rightarrow$)
 Since  ${}_RT \oplus * \cong {}_RR^t$ for some positive integer $t$,  we apply to this 
the functor $A \ot_R -$ from left $R$-modules into $A$-$B$-bimodules which results in 
${}_AA\! \ot_B\! A_B \oplus * \cong {}_AA_B^t$, after using the Galois $A$-$B$-isomorphism 
$A \ot_B A \cong  $ $A \ot_R T$. Hence, $A | B$ is right D2, and left D2 since $A|B$ is Frobenius \cite[6.4]{KS}.  

Let $\mathcal{E} := \End A_B$.  The module $A_B$ is balanced iff the natural bimodule ${}_{\mathcal{E}}A_B$
is faithfully balanced, which we proceed to show based on the following claim.
Let $R$ be a ring, $M_R$ and ${}_RV$ modules with ${}_RV$ finite projective.  If 
$\sum_j m_j \phi(v_j) = 0$ for all $\phi$ in the left $R$-dual ${}^*V$, then $\sum_j m_j \ot_R v_j = 0$.
This claim follows immediately by using dual bases $f_i \in {}^*V$, $w_i \in V$.  

Given $F \in \End {}_{\mathcal{E}}A$, it suffices to show that $F = \rho_b$ for some
$b \in B$.  Since $\lambda_a \in \mathcal{E}$, $F \circ \lambda_a = \lambda_a \circ F$ for all $a \in A$,
whence $F = \rho_{F(1)}$.  Designate $F(1) = a$.  If we show that $a\0 \ot a\1 = a \ot 1$ after applying
the right $T$-valued coaction on $A$,
then $a \in A^{\rm co \, T} = B$. For each $\alpha \in {}^*({}_RT)$, define $\overline{\alpha} \in \End A_B$
by $\overline{\alpha} (x) = x\0 \alpha(x\1)$. Since $\rho_r \in \mathcal{E}$ for each $r \in R$ by 2.4, 
$$a \alpha(1_T) = F(\overline{\alpha} (1)) = \overline{\alpha} (F(1)) = a\0 \alpha(a\1) $$
for all $\alpha \in {}^*T$. By the claim $a\0 \ot_R a\1 = a \ot 1_T$.

($\Leftarrow$) Let $T$ be the left projective right bialgebroid $(A \otimes_B A)^B$ over $R = C_A(B)$. 
Using a right D2 quasibase, we give $A$ the structure  of
 a right $T$-comodule algebra via $
 a\0 \ot a\1 := \sum_j \gamma_j(a) \ot u_j \in A \ot_R T$,
the details of which are in \cite[5.1]{K2004}. The D2 condition ensures that $\theta: A \ot_R T \to A \ot_B A$
defined by $\theta(a \ot t) = at^1 \ot t^2$ is an isomorphism.    

Note that for each $b\in B$
$$ b\0 \ot b\1 = \sum_j \gamma_j(b) \ot_R u_j = b \ot \sum_j \gamma_j(1)u_j = b \ot 1_T $$
so $B \subseteq A^{\rm co \, \rho}$.  
The converse: if $\rho(x) = x \ot 1_T$ $ = \sum_j \gamma_j(x) \ot u_j$ applying $\theta$ we
obtain $x \ot_B 1 = 1 \ot_B x$. Since
 $A_B$ is balanced, we know $A^S = B$ under the action $\lact$ of $S$ on $A$ \cite[4.1]{KS}.
Applying $\mu (\alpha \ot \id)$ for
each $\alpha \in S$, where $\mu$ is multiplication, we obtain $\alpha \lact x =$ $ \alpha(x) =$ $ \alpha(1) x$,
whence $x \in B$.
 
 The Galois  mapping $\beta: A \ot_B A \to A \ot_R T$
given by 
\begin{equation}
\label{eq: can}
\beta(a \ot a') := a {a'}\0 \ot_R {a'}\1
\end{equation}
 is an isomorphism  since $\theta$ is an inverse by eq.~(\ref{eq: right}).  \textit{Q.e.d.}

\begin{cor}
Let $k$ be a field and $A | B$ be a Frobenius extension of $k$-algebras 
with centralizer $R$
a separable $k$-algebra.  The extension $A | B$ is  weak Hopf-Galois
iff $A | B$ is D2 with $A_B$ balanced.  
\end{cor} 

The proof of $\Leftarrow$ depends first on recalling that the right $R$-bialgebroid
$T$ is a weak bialgebra since $R$ has an index-one Frobenius system
$(\phi: R \to k, e_i, f_i \in R)$ where $\sum_i e_i f_i = 1_A$
and $\sum_i \phi(re_i)f_i = r = \sum_i e_i \phi(f_i r)$ for all $r \in R$,
whence $\cop(t) = \sum_i t\1 e_i \ot_k f_i t\2$ and $\eps(t) = \phi(t^1 t^2)$
satisfy the axioms of a weak bialgebra \cite[(96)]{KS}.
Since $A | B$ is Frobenius, the dual bases tensor is a nondegenerate
right integral in $T$, whence $T$ is weak Hopf algebra by the Larson-Sweedler-Vecsernyes
theorem.  The coaction on $A$ has values
in $A \ot_R T \cong (A \ot_k T)\cop(1)$, given
by $a \ot t \mapsto \sum_i ae_i \ot f_i t$, an isomorphism
of the Galois $A$-corings  in \cite[5.1]{K2004}
and in \cite[2.1]{CG}. The proof of $\Rightarrow$ follows from
the fact that a weak Hopf-Galois extension is D2 
and an argument that $A_B$ is balanced like the one above.  

\begin{ex}
A separable field extension $\overline{k} | k$ is a weak Hopf-Galois extension, since $\overline{k}$ is
a separable $k$-algebra. 
\end{ex}

The theorem provides another proof and extends the theorems  \cite[8.14]{KS} and \cite[6.6]{KN} as we see below.
We define an \textit{irreducible} $k$-algebra extension to be an extension where the centralizer is the trivial $k1$.  
\begin{cor}  
Let $A | B$ be an irreducible extension.  
The extension $A | B$ is finite Hopf-Galois $\Longleftrightarrow$ $A | B$ is a D2, right balanced extension.  
\end{cor}
The proof of $\Rightarrow$ does not require the centralizer to be trivial.  The Frobenius condition may
be dropped here from the proof of $\Leftarrow$ since the bialgebra $T$ acts Galois implies it is
a Hopf algebra \cite{S}.
If the characteristic of $k$ is zero, the Larson-Radford theorem permits the condition 
``right balanced'' to be replaced by ``separable extension''
\cite[4.1]{K2004}.

\noindent 
Matematiska Institutionen, G{\" o}teborg 
University, 
S-412 96 G{\" o}teborg, Sweden, 
lkadison@c2i.net, 
September 20, 2004.  

\end{document}